\documentclass[10pt]{amsart}
\usepackage[a4paper, margin=1in, footskip=0.5in]{geometry}

\usepackage{a4, a4wide}
\usepackage{amssymb}
\usepackage{amsmath, bm}
\usepackage{amsthm}
\usepackage{amstext}
\usepackage{amscd}
\usepackage{stmaryrd}
\usepackage{latexsym}
\usepackage{mathtools}
\usepackage{graphics}
\usepackage{color}
\usepackage{enumitem}
\usepackage[all, cmtip]{xy}
\usepackage{enumitem}
\usepackage{xcolor}
\usepackage{chapterbib} 
\usepackage{footnote}
\usepackage{etoolbox}

\usepackage[colorlinks, pagebackref]{hyperref}
\hypersetup{
  colorlinks=true,
  citecolor=blue,
  linkcolor=blue,
  urlcolor=blue}
\usepackage{cleveref}

\usepackage{tikz-cd}
\tikzset{%
	symbol/.style={%
		draw=none,
		every to/.append style={%
			edge node={node [sloped, allow upside down, auto=false]{$#1$}}}
	}
}
 
\makeatletter
\def\@tocline#1#2#3#4#5#6#7{\relax
  \ifnum #1>\c@tocdepth 
  \else
    \par \addpenalty\@secpenalty\addvspace{#2}%
    \begingroup \hyphenpenalty\@M
    \@ifempty{#4}{%
      \@tempdima\csname r@tocindent\number#1\endcsname\relax
    }{%
      \@tempdima#4\relax
    }%
    \parindent\z@ \leftskip#3\relax \advance\leftskip\@tempdima\relax
    \rightskip\@pnumwidth plus4em \parfillskip-\@pnumwidth
    #5\leavevmode\hskip-\@tempdima
      \ifcase #1
      \or\or \hskip 2em \or \hskip 2homologyem \else \hskip 3em \fi%
      #6\nobreak\relax
    \dotfill\hbox to\@pnumwidth{\@tocpagenum{#7}}\par
    \nobreak
    \endgroup
  \fi}
\makeatother

\theoremstyle{plain}

\newtheorem{theorem}{Theorem}[section]

\newtheorem{lemma}[theorem]{Lemma}

\newtheorem{proposition}[theorem]{Proposition}
\newtheorem{conjecture}[theorem]{Conjecture}
\theoremstyle{definition}

\newtheorem{notations}[theorem]{Notations}
\newtheorem{remark}[theorem]{Remark}

\newtheorem{definition}[theorem]{Definition}
\newtheorem{example}[theorem]{Example}

\numberwithin{equation}{subsection}

\newcommand{\dlim}{\mathop{\varinjlim}\limits}  
\DeclareMathOperator{\holim}{holim} 

\newcommand{\Smk}{{\rm Sm/k}}
\newcommand{\SmS}{{\rm Sm/S}}

\newcommand{\Hom}{{\rm Hom}}

\newcommand{\im}{{\rm Im \ }}
\newcommand{\Spec}{{\rm Spec \,}}

\newcommand{\ie}{{\it i.e.\/},\ }

\newcommand{\nis}{Nis}
\newcommand{\zar}{Zar}
\newcommand{\et}{\rm \acute{e}t}

\newcommand{\N}{{\mathbb N}}

\newcommand{\A}{\mathbb A}

\usepackage[colorinlistoftodos, textsize=tiny]{todonotes}
\usepackage{mathrsfs}

\def\<{\langle}
\def\>{\rangle} 
\def\-{\overline} 
\def\~{\widetilde}
\def\^{\widehat}

\def\@{\mathcal}
\def\!{\mathscr}
\def\#{\mathbb}
\def\*{\mathbf}
\def\_{\underline}
 
\def\x{\times}

\def\.{\cdot}

\newcommand{\simps}{\rm{\Delta^{op}Shv}}

\input{xy}
\xyoption{all}

\begin{document}

\title{Non-finite type \'etale sites over fields}

\author{Sujeet Dhamore}
\address{Department of Mathematics, Indian Institute of Science Education and Research (IISER) Pune, Dr. Homi Bhabha road, Pashan, Pune 411008 India}
\email{sujeet.dhamore@students.iiserpune.ac.in}

\author{Amit Hogadi}
\address{Department of Mathematics, Indian Institute of Science Education and Research (IISER) Pune, Dr. Homi Bhabha road, Pashan, Pune 411008 India}
\email{amit@iiserpune.ac.in}

\author{Rakesh Pawar}
\address{UMPA, ENS de Lyon, UMR 5669, CNRS,
46 all\'ee d'Italie 69364 Lyon Cedex 07, France}
\address{Current Affiliation: Harish-Chandra Research Institute, A CI of Homi Bhabha National Institute, Chhatnag Road,
Jhunsi, Prayagraj 211019, India}
\email{rakesh.pawar@ens-lyon.fr, rakeshpawar@hri.res.in}

\date{\today}

\begin{abstract} We consider the notion of finite type-ness of a site introduced by Morel and Voevodsky, for the \'etale site of a field. 
For a given field $k$, we conjecture that the \'etale site of $\Smk$ is of finite type if and only if the field $k$ admits a finite extension of finite cohomological dimension. We prove this conjecture in some cases, e.g. in the case when $k$ is countable, or in the case when the $p$-cohomological dimension $cd_p(k)$ is infinite for infinitely many primes $p$. 
\end{abstract}

\maketitle


\section{Introduction}

In modern development of homotopy theory in algebraic geometry e.g., as in \cite{MV}, given a site with a Grothendieck topology, one works with the category of simplicial sheaves (or presheaves) on the site with local-injective model structure. 
In the category of simplicial sets, Kan's fibrant resolution functor exists such that it preserves finite limits and fibrations \cite[Lemma 1.67, page 70]{MV}. In contrast to this, in the category of simplicial sheaves (or presheaves) on a site, currently there is no known construction of such a fibrant resolution functor that preserves finite limits and fibrations. 

As noted in \cite[Theorem 1.66, pages 69-70]{MV}, if the site has enough points, the natural (Godement) functor  $$Ex^{\@G}: \simps(T)\to \simps(T)$$ is such a fibrant resolution functor under the additional assumption that the site $T$ is Postnikov complete/finite type. 
However, as remarked on \cite[page 66]{MV}, we do not know if the assumption of finite type can be avoided for the existence of such a  fibrant resolution functor.  

We recall notion of finite type site by Morel and Voevodsky.
\begin{definition}(Finite type site) \cite[Definition 1.31, page 58]{MV}\label{ft}
	A site $T$ is called of finite type, if for each simplicial sheaf $\@X$ on $T$, the canonical map $$\@X\to \holim_{n\geq0} Ex(P^n(\@X))$$ is a weak equivalence. Here, $P^n(\@X)$ is the $n^{th}$ Postnikov truncation (see \Cref{def:postnikov}) and $Ex$ is a resolution functor as part of the local-injective model structure.
\end{definition}

The Nisnevich site $\rm{(\SmS)_{\nis}}$, being of finite type (see~\Cref{rem:fin-type}), has various important consequences. For example,
\begin{enumerate}
\item  There exists an $\mathbb{A}^1$-fibrant replacement functor $L_{\mathbb{A}^1}$ on $\simps((Sm/S)_{\nis})$, which commutes with finite limits. This is a consequence of existence of a simplicial fibrant replacement functor $Ex^{\@G}$  which commutes with finite limits, and of ~\cite[Lemma 3.20, page 96]{MV}. In particular, given a simplicial group $\@G$, $L_{\#A^1}\@G$ is also a simplicial group and  given a $\@G$-torsor $\@E$, $L_{\#A^1}\@E$ is a  $L_{\#A^1}\@G$-torsor (\cite[Definition 6.49, Theorem 6.50]{morel}). For applications, reader can refer to \cite[Remark A.13]{morel} and in particular, it is used in the proofs of \cite[Theorem 6.46, Lemma B.7(2)]{morel}.

\item Maps in the homotopy category can be constructed inductively using obstruction theory (see~\cite[Appendix B]{morel}), particularly for studying vector bundles or principal bundles on a smooth affine scheme combined with representability theorems (as in \cite{AHW17}, \cite{AHW18}). 
\end{enumerate}

Note that while for the second point above, convergence of Postnikov tower is required only for the particular object in question, the first consequence requires the site to be of finite type. 
 
One would also like to have the above consequences in the \'etale topology.  In particular, one would like to know whether the category of smooth schemes over a field with  \'etale topology is of finite type. There are at least two possible reasons why we think the study of \'etale topology is interesting and potentially opens the doors to new applications which are not accessible in Nisnevich topology. 
 \begin{enumerate}
 \item  Let $X/k$ be a smooth projective variety. If $k$ is an uncountable field of characteristic zero and $X/k$ is rationally connected, then $X$ is \'etale-$\A^1$-connected.  
Conversely, under the assumption that the field $k$ is of characteristic 0, if a smooth projective variety $X$ is \'etale-$\A^1$-connected, then $X/k$ is rationally connected. For a detailed proof, see the Appendix. This provides a compelling reason to believe that \'etale version of $\A^1$-homotopy theory might be important for the study of rationally connected varieties at least over fields of characteristic zero. 
 \item The dichotomy between split and non-split algebraic groups appears to be less important in \'etale topology. To give an example,  in~\cite[Corollary 3.4]{BHS23}, it was shown that for a reductive algebraic group $G$, $Sing_*^{\#A^1}G$ is $\#A^1$-(Nis)-local iff $G$ is isotropic. In \'etale version of this theory, if $Sing_*^{\#A^1}G$ was known to be $\#A^1$-local for split groups, it would automatically be known for all algebraic groups, since every algebraic group is \'etale-locally split. 
 \end{enumerate}
 
This paper aims to make a modest contribution in the development of  \'etale version of $\#A^1$-homotopy theory by studying when the site $\rm{(Sm/k)_{\et}}$ of finite type smooth schemes over a field $k$ is of finite type. The goal of this paper is to give a criterion for this site to be of finite type. We conjecture the following.

\begin{conjecture}\label{conjecture}
	Let $k$ be any field. Let $\simps((Sm/k)_{\et})$ denote the category of \'etale simplicial sheaves on $(\Smk)_{\et}$ with local injective model structure. Then the site $\rm{(Sm/k)_{\et}}$ is of finite type if and only if there exists a finite extension $L/k$ such that $cd(L)<\infty$. 
\end{conjecture}
Here for a field $L$, and a prime $p$, the $p$-cohomological dimension is denoted by $cd_p(L)$ and the cohomological dimension is denoted by $cd(L)$.

The following is the main result of the paper which proves conjecture for countable fields. 

\begin{theorem}\label{ftcriterion+}
	Conjecture \ref{conjecture} holds, if the absolute Galois group $G_k = Gal(k^{sep}/k)$ is first-countable. In particular, if the field $k$ is countable, then $G_k$ is first-countable, hence conjecture \ref{conjecture} holds. 
\end{theorem}

For any field k, the \textbf{if} part of the \Cref{conjecture} follows from the following result of Morel and
Voevodsky, which gives a sufficient cohomological criterion for a site to be of finite type.

\begin{theorem}\cite[Theorem 1.37, page 60]{MV}\label{fin-type}
	Let $T$ be a site and suppose that there exists a family $(A_d)_{d \geq 0}$ of classes of objects of $T$ such that the following conditions hold:
	\begin{enumerate}
		\item Any object $U$ in $A_d$ has cohomological dimension $\leq d$ i.e. for any sheaf $F$ on $T/U$ and any $m>d$ one has $H^m(U,F)=0$.
		\item For any object $V$ of $T$ there exists an integer $d_V$ such that any covering of $V$ in $T$ has a refinement  of the form $\{U_i\to V\}$ with $U_i$ being in $A_{d_V}$.
	\end{enumerate}
	Then $T$ is a site of finite type.
\end{theorem}

\begin{remark}\label{rem:fin-type}
	As a consequence of Theorem \ref{fin-type}, for a base scheme $S$ of finite Krull dimension, the sites $\rm{(\SmS)_{\zar}}$ and $\rm{(\SmS)_{\nis}}$ are of finite type (see \cite[Proposition 1.8, page 98]{MV}).
\end{remark}

Note that the Theorem \ref{fin-type} implies that, if $k$ is a field, such that there exists a finite extension of $k$ which is of finite cohomological dimension, then $\rm{(Sm/k)_{\et}}$ is of finite type.
Hence the \textbf{only if} part of the Conjecture \ref{conjecture} is open.

We also show that, certain cohomological boundedness conditions on the field $k$ are necessary for the site $\rm{(Sm/k)_{\et}}$ to be of finite type. In following results, we show \textbf{only if} part of Conjecture \ref{conjecture} for some cases.

\begin{theorem}\label{thmlargeorder}
Let $k$ be any field which satisfies the condition that for each $n\in\N$, there exists $G_k$-module $A_n$ such that 
$$\limsup_{n\to \infty} exp(H^n(G_{k},A_n))=\infty$$
where for an abelian group $B$,  the exponent of $B$ is denoted by  $exp(B)$. Then the site $\rm{(Sm/k)_{\et}}$ is not of finite type.
\end{theorem}

Observe that, a field $k$ which either has $cd_p(k)=\infty$ for infinitely many primes $p$ or has arbitrarily large $cd_p(k)$ for arbitrarily large primes $p$ satisfies hypothesis of Theorem \ref{thmlargeorder}. This is summarized in the following result.

\begin{proposition}\label{cdpcondition}
Let $k$ be a field such that $\limsup cd_p(k) = \infty$ where limit is over all primes $p$. Then the site $\rm{(Sm/k)_{\et}}$ is not of finite type.
\end{proposition}

Note that, if a field $k$ satisfies assumptions of either Proposition \ref{cdpcondition} or Theorem \ref{thmlargeorder} then for each  finite extension $L/k$,  $cd(L)=\infty$.

\subsection*{Acknowledgments}
The first author was supported by Ph.D. fellowship (File No: 09/936(0247)/2019-EMR-I) of Council of Scientific \& Industrial Research (CSIR), India. The third author would like to thank Matthias Wendt for helpful discussions. The third author was supported by the French ANR
project: Motivic homotopy, quadratic invariants and diagonal classes (ANR-21-CE40-0015).


\section{preliminaries}

In this section, we fix some notations and recall definitions from \cite{MV}.

For any Grothendieck site $T$ with enough points, $\simps(T)$ and $\simps_{Ab}(T)$ denote the category of simplicial sheaves and the category of simplicial abelian sheaves on the site $T$ respectively.  We will consider $\simps(T)$ and $\simps_{Ab}(T)$ as a simplicial model categories with the local injective model structure as endowed in \cite[pages 48-49]{MV}.

Analogues to the situation for simplicial sets, Morel and Voevodsky introduced the notion of (simplicial) Postnikov truncations of a given simplicial sheaf.
\begin{definition}\cite[Page 57]{MV}\label{def:postnikov}
	Given a simplicial sheaf $\@X \in \simps(T)$ for each $n\geq 0$, $P^n(\@X)$ is the sheafification of the presheaf given by 
	\begin{center}
		$U \longmapsto P^{(n)}(\@X(U)) = \im(\@X(U)\longrightarrow cosk_n(\@X(U)))$ for each $U \in T$
	\end{center}
	where $P^{(n)}(S)$ denotes $n^{th}$ Moore-Postnikov truncation for any simplicial set $S$.
\end{definition}

\begin{notations} We use the following notations throughout the paper.
	\begin{enumerate} 
		
		\item Given a simplicial sheaf $\@X$, the homotopy sheaf $\pi_0^s(\@X)$ is defined as the sheafification of the presheaf $\pi_0^{pre}(\@X)$ given as 
		\begin{center}
			$\pi_0^{pre}(\@X)(U) = \pi_0(\@X(U))$ for each $U \in T$
		\end{center}
		where $\pi_0(S)$ denotes the set of connected components of any simplicial set $S$.
		
		\item Another description of $\pi_0^{pre}(\@X)$ is given as $U \longmapsto \Hom_{\@H_s}(U,\@X)$, where $\@H_s$ is the (simplicial)-homotopy category of $\simps(T)$.
		
		\item For a given site $T$, let $Ex: \simps(T) \to \simps(T)$ denote a resolution functor or a fibrant replacement functor as defined in \cite[Definition 1.6, page 49]{MV}.
		
		\item Given a field $k$ and a fixed separable closure $k^{sep}/k$, for any intermediate field extension $k^{sep}/L/k$, let $G_L := Gal(k^{sep}/L)$ be absolute Galois group of $L$. 
		
	\end{enumerate}
\end{notations}


\section{Main Results}
The goal of this section is to prove Theorems \ref{ftcriterion+}, \ref{thmlargeorder} and Proposition \ref{cdpcondition}. Lemma \ref{ftcondition} is at the heart of the proofs of these results. The broad idea of the proof of this lemma, is borrowed from the following example of Morel and Voevodsky.
\begin{example}(\cite[Example 1.30, page 58]{MV}) 
For $G=\prod_{i>0} \#Z/2$, the site $(G\hbox{-}sets)_{df}$ of finite $G$-sets (with jointly surjective coverings) is not of finite type. This relies on the fact that $H^i(\#Z/2, \#Z/2)\neq 0$ for all $i>0$, or in other words, $\#Z/2$ has infinite (group) cohomological dimension. 
\end{example}
 
\begin{lemma} \label{ftcondition} Let $k$ be a field. Assume that for each $i\in \N$, there is a continuous $G_k$-module $A_i$ and classes $$\{\alpha_i \in H^i(G_k, A_i)\}_{i\in \N}$$ such that for each $L\supset k$ finite Galois extension, there is $m$ (possibly depending on $L$) such that $(\alpha_{m})_{|L}\neq 0$ in $H^{m}(G_L,A_{m})$. Then the site $\rm{(Sm/k)_{\et}}$ is not of finite type.
\end{lemma}
\begin{proof} To prove the lemma we need to construct an \'etale simplicial sheaf $\@X$ such that $\@X\to \holim_{n\geq0} Ex(P^n(\@X))$ is not a weak equivalence. We do this in several steps. 
\\
\_{\it Step 1:} In this step we construct the required simplicial sheaf $\@X$.
 
 For each $G_k$-module $A_i$ we have the associated abelian sheaf $\~{\@A_i}\in{\rm Shv_{Ab}((\Spec k)_{\et})}$. The inclusion functor $f:(\Spec k)_{\et} \rightarrow \rm{(Sm/k)_{\et}}$ induces the adjoint pair between respective categories of sheaves by \cite[Chapter I, 3.6]{tamme} denoted as,
$$f^{-1}: {\rm Shv_{Ab}((\Spec k)_{\et})}\rightleftarrows {\rm Shv_{Ab}((Sm/k)_{\et})}: f_*.$$

Let $\@A_i:=f^{-1}(\~{\@A_i})$ be the sheaves on $\rm{(Sm/k)_{\et}}$.

We recall classical Dold-Kan correspondence, \ie an equivalence between the category ${\rm sAb}$ of simplicial abelian groups  and the category  ${\rm Ch_{\geq 0}}$ of  chain complexes (in non-negative degrees) of abelian groups (see \cite[Chapter 2, 2.3]{GJ})
$$ N: {\rm sAb} \rightleftarrows {\rm Ch_{\geq 0}}: \Gamma.$$
By \cite[Proposition 1.24, page 56]{MV}, for site $T$, by applying the functors $N$ and $\Gamma$ sectionwise, this equivalence extends to the equivalence between the category ${\rm Shv_{sAb}(T)}$ of sheaves of simplicial abelian groups  and the category ${\rm Shv_{Ch_{\geq 0}}(T)}$ of sheaves of complexes (in non-negative degrees) of abelian groups 
$$ \~N: {\rm Shv_{sAb}(T)} \rightleftarrows {\rm Shv_{Ch_{\geq 0}}(T)}: \~\Gamma .$$
More precisely, for $\@F \in {\rm Shv_{sAb}(T)}$, the sheaf $\~N\@F$ is given as $U \mapsto N(\@F(U))$ for all $U \in T$, and for $\@G \in {\rm Shv_{Ch_{\geq 0}}(T)}$, the sheaf $\~\Gamma\@G$ is given as $V \mapsto \Gamma(\@G(V))$ for all $V \in T$.

Consequently, we obtain a model of Eilenberg-MacLane sheaves corresponding to the complex $\@A_i[i]$ with $\@A_i$ in degree $i$, defined as $\@K(\@A_i,i) = \~\Gamma(\@A_i[i])\in \simps((Sm/k)_{\et})$, which is described as
\begin{center}
	$ \@K(\@A_i,i)(U) = \~\Gamma(\@A_i[i])(U) = \Gamma(\@A_i(U)[i]) = K(\@A_i(U),i)$ for any $U \in \rm{(Sm/k)_{\et}}.$ 
\end{center}
Note that $K(\@A_i(U),i)$ is the Eilenberg-MacLane simplicial abelian group corresponding to the complex of abelian groups $\@A_i(U)[i]$. We define $\@X:= \prod_{i>0}\@K(\@A_i,i) \in \simps((Sm/k)_{\et})$ as the product of simplicial sheaves. \\~\\
	\_{\it Step 2:} In this step we show that the simplicial sheaf $\@X$ is connected.  
	
For each $U \in \rm{(Sm/k)_{\et}}$, the homotopy presheaf $\pi_0^{pre}(\@X)$ is given as
\begin{align*}
	\pi_0^{pre}(\@X)(U) = \pi_0(\@X(U)) &= \pi_0(\prod_{i>0}\@{K}(\@A_i,i)(U)) \\
	&= \pi_0(\prod_{i>0}K(\@A_i(U),i)) \\
	&= \prod_{i>0} \pi_0(K(\@A_i(U),i)).
\end{align*}
For the last equality, first we observe that each $K(\@A_i(U),i)$ is a simplicial abelian group hence a Kan complex \cite[Chapter 1, Lemma 3.4]{GJ}. Also, the functor $\pi_0(-)$ on simplicial sets commutes with infinite product of Kan complexes. To see this, consider Kan complexes $K_i$ for $i \in I$, then clearly $\prod_{i\in I}K_i$ is also a Kan complex. The natural map $\pi_0(\prod_{i\in I}K_i) \xrightarrow{f} \prod_{i\in I} \pi_0(K_i)$ is surjective, since product of 0-simplices (representing connected components in $K_i$) represents a connected component of  $\prod_{i\in I}K_i$. For injectivity, if two 0-simplices $x=(x_i)$ and $y=(y_i)$ of $\prod_{i\in I}K_i$ are such that, all $x_i$ and $y_i$ are in the same connected component of $K_i$, then there is 1-simplex $s_i$ in each $K_i$ which has faces as $x_i$ and $y_i$ (since each $K_i$ is Kan complex). Then the tuple $(s_i)_{i \in I}$ is a 1-simplex of $\prod_{i\in I}K_i$ implying that $x$ and $y$ are in the same connected component of $\prod_{i\in I}K_i$.

Since $K(\@A_i(U),i)$ is connected (for $i>0$), we have that for each $U$, $\pi_0^{pre}(\@X)(U) = *$. Hence, the homotopy presheaf $\pi_0^{pre}(\@X)=*$, which implies that the associated sheaf $\pi_0^s(\@X)=*$. \\~\\
	\_{\it Step 3:} We determine Postnikov truncations $P^n(\@X)$ corresponding to $\@X$. 
	
	Recall from \Cref{def:postnikov} that $P^n(\@X)$ is the sheaf associated to the presheaf 
$$U \mapsto P^{(n)}(\@X(U))=\im(\@X(U)\xrightarrow{\eta} cosk_n(\@X(U)))$$
where $\eta$ for any simplicial set $X$ is the canonical map from $X$ to its coskeleton $cosk_n(X)$ such that for any $m$-simplex $x$ given as a map $\Delta^m \to X$
\begin{center}
$\eta(x) = x_{|sk_n(\Delta^m)}$ given as the map $sk_n(\Delta^m) \to X$.
\end{center}
Observe that, since $cosk_n$ is right adjoint \cite[Example 1.1.9]{Riehl}, it commutes with limits and we have 
$$cosk_n(\@X(U))=cosk_n(\prod_{i>0}\@{K}(\@A_i,i)(U))=\prod_{i>0}cosk_n(K(\@A_i(U), i)).$$ 
Hence, $\eta$ is given as 
$$ \prod_{i>0} K(\@A_i(U), i) \xrightarrow{\eta = \prod_i{\eta_i}} \prod_{i>0}cosk_n(K(\@A_i(U), i)) $$
where $\eta_i : K(\@A_i(U), i) \to cosk_n(K(\@A_i(U), i))$ is the canonical map to coskeleton. \\
Also observe that $cosk_n(K(\@A_i(U), i))=*$ for $i>n$. Hence we have
$$P^{(n)}(\@X(U)) = \prod_{i=1}^{n}\im(K(\@A_i(U), i) \to cosk_n(K(\@A_i(U), i))) = \prod_{i=1}^{n} P^{(n)}(\@{K}(\@A_i(U),i)).$$
Since $\@{K}(\@A_i(U),i)$ is weakly equivalent to $P^{(n)}(\@{K}(\@A_i(U),i))$ for $i\leq n$, we get 
$\prod_{i=1}^n\@{K}(\@A_i,i)(U)$ is weakly equivalent to $P^{(n)}(\@X(U))$ compatible with $n$. Hence, there is weak equivalence of simplicial sheaves $\prod_{i=1}^{n} \@K(\@A_i, i) \rightarrow P^n(\@X)$. 
\\~\\
\_{\it Step 4:} We show that simplicial sheaf $\prod_{i>0} Ex(\@K(\@A_i, i))$ is weakly equivalent to $\holim_n Ex(P^n(\@X))$. To show that the site $\rm{(Sm/k)_{\et}}$ is not of finite type, it is enough to show that the homotopy sheaf $\pi_0^s(\holim_n Ex(P^n(\@X))) \cong \pi_0^s(\prod_{i>0} Ex(\@K(\@A_i, i)))$ is non trivial. 
	
Observe that the map $ Ex(\prod_{i=1}^{n} \@K(\@A_i, i)) \to Ex(P^n(\@X))$ is sectionwise weak equivalence of simplicial sets by \cite[Corollary 5.13]{jardine}. The homotopy limit of simplicial sheaves is defined sectionwise, then by Homotopy lemma \cite[XI, 5.6]{BK} the induced map
$$\holim_n Ex(\prod_{i=1}^{n} \@K(\@A_i, i)) \longrightarrow \holim_n Ex(P^n(\@X))$$ is weak equivalence of simplicial sheaves.
Further, we have $$Ex(\prod_{i=1}^{n} \@K(\@A_i, i)) \longleftarrow \prod_{i=1}^{n} \@K(\@A_i, i) \longrightarrow \prod_{i=1}^{n} Ex(\@K(\@A_i,i))$$ where both maps are weak equivalences, first is part of definition of resolution functor $Ex$ and second is obtained as finite product of weak equivalences $\@K(\@A_i,i) \rightarrow Ex(\@K(\@A_i,i))$.
Then using lifting properties of model structure, for each $n>0$, we can induce weak equivalences
$$ Ex(\prod_{i=1}^{n} \@K(\@A_i, i)) \xrightarrow{f_n} \prod_{i=1}^{n} Ex(\@K(\@A_i,i)) $$ which are compatible with $n$. Then by \cite[XI, 5.6]{BK} we get weak equivalence
 $$\holim_n Ex(\prod_{i=1}^{n} \@K(\@A_i, i)) \longrightarrow \holim_n \prod_{i=1}^{n} Ex(\@K(\@A_i,i)).$$
 
Observe that the sequence $\big\{\prod_{i=1}^{n} Ex(\@K(\@A_i,i))\big\}_n$ is a tower of fibrant objects and maps are projections which are fibrations. Then using \cite[XI, 4.1]{BK} the natural map $$\lim_n \prod_{i=1}^{n} Ex(\@K(\@A_i,i)) \longrightarrow \holim_n \prod_{i=1}^{n} Ex(\@K(\@A_i,i))$$ is a weak equivalence.
Let us denote $\~{\@X}:=\prod_{i>0} Ex(\@K(\@A_i, i)) = \lim_n \prod_{i=1}^{n} Ex(\@K(\@A_i,i))$. \\~\\
	\_{\it Step 5:} We show that the homotopy sheaf $\pi_0^s(\~{\@X}) \neq *$. 
	
For each $U\in \rm{(Sm/k)_{\et}}$, the homotopy presheaf is given as
$$\pi_0^{pre}(\~{\@X})(U) = \Hom_{\@H_s}(U, \prod_{i>0} Ex(\@K(\@A_i, i))) =\prod_{i>0} H^i(U, \@A_i).$$
where last equality follows from \cite[Proposition 1.26, page 57]{MV}.
Then the stalk of $\pi_0^s(\~{\@X})$ at geometric point $\bar{s} = \Spec (k^{sep})$ is given by
\begin{align*} 
	\pi_0^s(\~{\@X})_{\bar{s}} = \pi_0^{pre}(\~{\@X})_{\bar{s}} &= \dlim_{U} \pi_0^{pre}(\~{\@X})(U) \\
		&= \dlim_{U} \prod_{i>0} H^i(U, \@A_i) \\
		&= \dlim_{k \subset L} (\prod_{i>0} H^i(\Spec L, \@A_i))
\end{align*}	
where the colimit runs over \'etale neighborhoods $U$ of $\bar{s}$ which are given by finite Galois extensions $L/k$ so that $k\subset L\subset k^{sep}$. Observe that,
\begin{center}
$H^i(\Spec L, \@A_i) \cong H^i(\Spec L, \~{\@A_i}) \cong H^i(G_L,A_i)$ for all $i>0$
\end{center}
where first isomorphism is by \cite[Chatper I, Corollary 3.9.3]{tamme} and second isomorphism is by our choice of $\~{\@A_i}$. Hence we have
$$\pi_0^s(\~{\@X})_{\bar{s}} \cong \dlim_{k \subset L} (\prod_{i>0} H^i(G_L, A_i)). $$
Let $\alpha:=\prod_i \alpha_i \in \prod_{i>0 }H^i(G_k, A_i)$. Then $\alpha|_L \neq 0 \ \text{in} \ \prod_{i>0} H^i(G_L, A_i)$ for each finite Galois extension $L/k$. Thus, $\alpha$ survives in above colimit, hence $\pi_0^s(\~{\@X})_{\bar{s}}\neq *$. This implies that the homotopy sheaf $\pi_0^s(\~{\@X})$ is nontrivial. 

Thus, the site $\rm{(Sm/k)_{\et}}$ is not of finite type.
\end{proof} 

\begin{remark}
The hypothesis of Lemma \ref{ftcondition} also implies that small \'etale site $(\Spec k)_{\et}$ is not of finite type.
\end{remark}

\begin{proof}[\textbf{Proof of \Cref{ftcriterion+}}]
By the assumption, the absolute Galois group $G_k$ is first-countable, so we can choose countable nested local basis at identity as open normal subgroups $\{G_{L_i}\}_{i \in \N}$ which correspond to sequence of finite Galois extensions  
$$ k \subset L_1 \subset L_2 \subset \dots \subset L_n \subset \dots \subset k^{sep} $$
such that for any $E/k$, a finite Galois extension, $E \subset L_m$ for some $m$.

We assume that for any finite extension $L/k$, we have $cd(L)=\infty$. For each $n \in \N$, $cd(L_n)=\infty$, so we can choose non-zero class $\alpha_n \in H^n(G_{L_n},A_n)$ for some $G_{L_n}$-module $A_n$. 

Let $H$ be a closed subgroup of a profinite group $G$ and let $A$ be a discrete $H$-module. Let $\hbox{Coind}^G_H(A)$ denote the Coinduced $G$-module (as in~\cite[page 242, Section 6.10]{RZ}).

By Shapiro's lemma \cite[Theorem 6.10.5]{RZ}, for $n \geq 0$ we have an isomorphism 
$$ H^n(G,\hbox{Coind}^G_H(A)) \cong H^n(H,A)$$
which factors through the restriction map \ie for each $n\geq 0$, the following diagram commutes
\[ \xymatrix{
	H^n(G,\hbox{Coind}^G_H(A)) \ar[rr]^{Shapiro} \ar[rd]^{res} & & H^n(H,A) \\
	& H^n(H, \hbox{Coind}^G_H(A)) \ar[ru]^{\nu} &
} \]
where $\nu$ is induced by the counit map $\hbox{Coind}^G_H(A) \rightarrow A$. Hence, the restriction map above is injective.

In particular, we have Coinduced $G_k$-modules, $M_n(A_n) := \hbox{Coind}^{G_k}_{G_{L_n}}(A_n)$. For each $\alpha_n$, we denote corresponding class through Shapiro's isomorphism in $H^n(G_k,M_n(A_n))$ by $\alpha_n$.

Since the restriction map above is injective, for each $n\in\N$, $(\alpha_n)_{|L_n} \neq 0$ in $H^n(G_{L_n},M_n(A_n))$. Furthermore, given a finite Galois extension  $E/k$, there is $m$ such that $E \subset L_m$. Under the composition 
$$H^m(G_k, M_m(A_m)) \to H^m(G_E, M_m(A_m)|_E)\to H^m(G_{L_m},M_m(A_m))$$
the element $\alpha_m\mapsto (\alpha_m)_{|L_m}\neq 0$, hence the image of $\alpha_m$ under $$H^n(G_k, M_n(A_n)) \to H^n(G_E, M_n(A_n)|_E)$$ is non-zero. Thus, the collection $\{\alpha_n\in  H^n(G_k,M_n(A_n))\}_n$ satisfies the hypothesis of Lemma \ref{ftcondition}, hence the site $\rm{(Sm/k)_{\et}}$ is not of finite type. 

Converse follows from Theorem \ref{fin-type}. 
\end{proof}

\begin{proof}[\textbf{Proof of \Cref{thmlargeorder}}]
	According to the hypothesis, there is a subsequence $\{n_i\}_{i\in\N}$ such that 
	$$\lim_i exp(H^{n_i}(G_k,A_{n_i})) = \lim_i c_i = \infty$$
	where sequence $\{c_i\}_{i \in \N}$ is either consists of $\infty's$ only or strictly monotonically increasing. We choose cohomology classes $\alpha_{n_i} \in H^{n_i}(G_k,A_{n_i})$ such that $|\alpha_{n_{i+1}}|>|\alpha_{n_i}|$ and $\lim_i |\alpha_{n_i}| = \infty$ (where for an abelian group $B$ and element $\beta \in B$, the order of $\beta$ is denoted by $|\beta|$ )
	
	Consider any finite Galois extension $L/k$, then the index $[G_k:G_L]=m$ for some $m\in\N$. Consider $n_i$ such that $|\alpha_{n_i}|>m$. Under the composition $$H^{n_i}(G_k, A_i) \xrightarrow{res} H^{n_i}(G_L, A_i) \xrightarrow{cor} H^{n_i}(G_k,A_i) $$  the element $\alpha_{n_i} \mapsto m\alpha_{n_i} \neq 0$, hence $(\alpha_{n_i})_{|L} \neq 0$.
	
		Thus, we have a collection $\{\alpha_n \in H^n(G_k, A_n)\}_{n\in \N}$, possibly with $\alpha_n = 0$ and $A_n = 0$ whenever $n \neq n_i$ for any $i$, which satisfies hypothesis of Lemma \ref{ftcondition}, hence the site $\rm{(Sm/k)_{\et}}$ is not of finite type.
\end{proof}

\begin{proof}[\textbf{Proof of \Cref{cdpcondition}}]
	By the hypothesis, there is a subsequence of primes $\{p_i\}_{i \in \N}$ such that $\lim_i cd_{p_i}(k)=\infty$ and sequence \{$cd_{p_i}(k)\}_{i \in \N}$ either consists of $\infty's$ only or is strictly monotonically increasing. 
	
	For a prime $p_1$ there is a $G_k$-module $A_1$ and a cohomology class $\alpha_{n_1}\in H^{n_1}(G_k,A_1)$ such that $|\alpha_{n_1}| \geq p_1$. For a prime $p_2$, since $cd_{p_2}(k)=\infty$ or $cd_{p_2}(k)>cd_{p_1}(k)$, we can choose $G_k$-module $A_2$ and class $\alpha_{n_2}\in H^{n_2}(G_k,A_2)$ such that $|\alpha_{n_2}| \geq p_2$ and $n_2 > n_1$. Inductively on $i$, we have classes $\{\alpha_{n_i} \in H^{n_i}(G_k, A_i)\}_{i\in \N}$ such that $|\alpha_{n_i}| \geq p_i$ and $n_{i+1} > n_i$ for each $i$. Thus the hypothesis of \Cref{thmlargeorder} is satisfied and the site $\rm{(Sm/k)_{\et}}$ is not of finite type.
\end{proof}

\section{Appendix: Rational connectedness versus \'etale-$\#A^1$-connectedness}\label{Appx}
The goal of this appendix is to  elaborate on the relationship between a rationally connected variety  over a field of characteristic zero and \'etale-$\A^1$-connectedness for a smooth projective variety. We give a detailed proof due to the lack of the reference in the literature. We refer the reader to \cite{carolina_kollar} and~\cite{kollar96} for general results on very free rational curves, space of morphisms and smoothness of the space of morphisms.
\begin{theorem} Let $X/k$ be a smooth projective variety. If $k$ is an uncountable field of characteristic zero and $X/k$ is rationally connected, then $X$ is \'etale-$\A^1$-connected.  
Conversely, under the assumption that the field $k$ is of characteristic 0, if a smooth projective variety $X$ is \'etale-$\A^1$-connected, then $X/k$ is rationally connected.  
\end{theorem} 
\begin{proof}
It is enough to show that for any strictly Henselian local scheme $U$ and maps $f,g:U\to X$, there is a map $H:U\times {\mathbb P}^1 \to X$ such that $H_0=f$ and $H_1=g$. Let $u\in U$ be the closed point which is of the form $\Spec E$ where $E$ is an algebraically closed field. Since $X$ is rationally connected and $E$ is algebraically closed, by~\cite[Page 37, Definition-Theorem 29]{carolina_kollar} there is a very free rational curve $h:{\mathbb P}^1_E\to X_E$ such that $h(0)= f(u)$ and $h(1)=g(u)$ \ie $$h^*(T_{X_E})\simeq \bigoplus_{i=0}^{\dim X_E} \@O_{\#P_E^1}(a_i)$$ where $a_i\geq 1$ for all $i$. Hence $H^1(\#P^1_E, h^*(T_{X_E})(-2))=0$ (by~\cite[Theorem 5.1]{Hart77}).

Let $B:=U\x \{0, 1\}\subset \#P^1_U$ and $\psi: B\to X_U$ be the $U$-morphism defined by $f$ on $U\x 0$ and the morphism $g$ on $U\x 1$. Consider the $U$-scheme of morphisms                         
$\Hom({\mathbb P}^1, X, \psi)$ defined by  the functor
$$ \Hom(\#P^1, X, \psi)(T):=\{ T\text{-morphism} \ h: T\times_U {\mathbb P}^1_U \to T\times_U X_U \ \text{such that} \ h|_{ T\x_U B}=\psi|_{T\x _U B} \}$$
for $U$-scheme $T$ (see \cite[chapter II, section II.1, 1.4, 1.5, 1.6]{kollar96}).  By above discussion the map $h: {\mathbb P}^1_E\to X_E$ corresponds to $[h]\in  \Hom(\#P^1, X, \psi)(E)$. By~\cite[Chapter II, Theorem 1.7]{kollar96} since $H^1(\#P^1_E, h^*(T_{X_E})\otimes I_{B_E})\simeq H^1(\#P^1_E, h^*(T_{X_E})(-2))=0$, the map $$ \Hom(\#P^1, X, \psi)\to U$$ is smooth at the $E$-rational point $[h]$. Since $U$ is strictly Henselian, there is a lift $H: U\to  \Hom(\#P^1, X, \psi)$ such that $H|_{\Spec E}: \Spec E\to U\to \Hom(\#P^1, X, \psi)$ is $h$. Thus, $$H :U\to  \Hom(\#P^1, X, \psi)$$ is such that $H|_{U\x 0} : U\x\{0\} \to U\x\#P^1\to U\x X$ is equal to $f$ and also   $H|_{U\x 1} : U\x\{1\} \to U\x\#P^1\to U\x X$ is equal to $g$.

Converse follows from~\cite[Definition 2.2.2 and Corollary 2.4.4]{AM11}. If $X$ is \'etale-$\A^1$-connected, then it is $\#A^1$-chain-connected over an algebraically closed field. Hence $X$ is rationally chain connected and over the field $k$ of characteristic 0, by~\cite[Corollary 4.28]{Deb01} it is rationally connected.

\end{proof}



\end{document}